\documentclass[12pt]{amsart}

\usepackage{mathabx}
\input bookman.sty
\boldmath

\usepackage{xcolor}

\usepackage{helvet}

\usepackage{graphics}
\usepackage{amssymb}
\usepackage{amsxtra}
\usepackage{amsmath}
\usepackage{mathrsfs}

\usepackage[arrow, matrix]{xy}
\xyoption{frame}

\theoremstyle{plain}

\newtheorem{theo}{Theorem}[section]
\newtheorem{lemm}[theo]{Lemma}
\newtheorem{prop}[theo]{Proposition}
\newtheorem{coro}[theo]{Corollary}

\theoremstyle{definition}

\newtheorem{defi}[theo]{Definition}
\newtheorem{rema}[theo]{Remark}

\newfont{\rmm}{cmr10 scaled 1000}

\newfont{\itt}{cmsl10 scaled 1000}

\newfont{\rM}{cmr10 scaled 1700}

\newcounter{lemma}[section]

\newcounter{tempcounter}

\newcommand{\lb}{\label}

\newcommand{\rrf}[1]{(\ref{#1})}

\renewcommand{\a}{\alpha}
\renewcommand{\b}{\beta}

\renewcommand{\d}{\delta}
\newcommand{\e}{\epsilon}

\renewcommand{\t}{\theta}
\renewcommand{\l}{\lambda}

\renewcommand{\r}{\rho}

\newcommand{\s}{\sigma}

\renewcommand{\L}{\Lambda}

\newcommand{\Si}{\Sigma}

\newcommand{\MM}{{\mathcal M}}
\newcommand{\NN}{{\mathcal N}}

\newcommand{\nn}{{\mathbb{N}}}

\newcommand{\qq}{{\mathbb{Q}}}
\newcommand{\rr}{{\mathbb{R}}}

\newcommand{\zz}{{\mathbb{Z}}}

\newcommand{\RRR}{{\mathbf{R}}}

\newcommand{\kkrest}{\begin{picture}(14,14)
\put(00,04){\line(1,0){14}}
\put(00,02){\line(1,0){14}}
\put(06,-4){\line(0,1){14}}
\put(08,-4){\line(0,1){14}}
\end{picture}     }

\newcommand{\krest}{~\kkrest~}

\newcommand{\Wh}{\text{\rm Wh}}
\newcommand{\Ker}{\text{\rm Ker }}

\newcommand{\Hom}{\text{\rm Hom}}

\renewcommand{\Im}{\text{\rm Im }}

\newcommand{\Int}{\text{\rm Int }}

\newcommand{\Id}{\text{\rm Id}}

\newcommand{\bere}{\begin{rema}}
\newcommand{\bede}{\begin{defi}}

\renewcommand{\beth}{\begin{theo}}
\newcommand{\bele}{\begin{lemm}}
\newcommand{\bepr}{\begin{prop}}
\newcommand{\beeq}{\begin{equation}}
\newcommand{\bega}{\begin{gather}}
\newcommand{\begaa}{\begin{gather*}}
\newcommand{\been}{\begin{enumerate}}

\newcommand{\bedee}{\begin{defii}}
\newcommand{\bethh}{\begin{theoo}}
\newcommand{\belee}{\begin{lemmm}}
\newcommand{\beprr}{\begin{propp}}

\newcommand{\beco}{\begin{coro}}

\newcommand{\beal}{\begin{aligned}}

\newcommand{\enre}{\end{rema}}

\newcommand{\enco}{\end{coro}}
\newcommand{\enpr}{\end{prop}}
\newcommand{\enth}{\end{theo}}
\newcommand{\enle}{\end{lemm}}
\newcommand{\enen}{\end{enumerate}}
\newcommand{\enga}{\end{gather}}
\newcommand{\engaa}{\end{gather*}}
\newcommand{\eneq}{\end{equation}}
\newcommand{\enal}{\end{aligned}}

\newcommand{\bq}{\begin{equation}}
\newcommand{\bqq}{\begin{equation*}}

\renewcommand{\leq}{\leqslant}
\renewcommand{\geq}{\geqslant}

\newcommand{\wi}{\widetilde}

\newcommand{\ove}{\overline}

\newcommand{\wh}{\widehat}

\newcommand{\sm}{\setminus}

\newcommand{\sbs}{\subset}

\newcommand{\tens}[1]{\underset{#1}{\otimes}}

\newcommand{\Lxi}{{\wh \L}_\xi}

\newcommand{\wrt}{with respect to}
\newcommand{\ho}{homomorphism}

\newcommand{\nei}{neighbourhood}

\newcommand{\Prf}{{\it Proof.\quad}}

\newcommand{\smo}{C^{\infty}}

\newcommand{\chart}{\Phi_p:U_p\to B^n(0,r_p)}
\newcommand{\atlas}{\{\Phi_p:U_p\to B^n(0,r_p)\}_{p\in S(f)}}

\newcommand{\pr}{\partial}

\newcommand{\qs}{\hfill\square}

\newcommand{\pa}{\vskip0.1in}

\newcommand{\mf}[2]{\mi {\phi_0}{\b_{#1}}{\b_{#2}}}

\newcommand{\arrh}[3]
{
\xymatrix{
{#1} \ar[r]^<<<<{#2}  &{#3}
}
}

\newcommand{\arrr}[1]
{\arrh {}{#1}{}}

\newcommand{\arrto}
{\xymatrix{{} \ar@{|-{>}}[r]  & {} } }

\newcommand{\arrinto}
{\xymatrix{{} \ar@{^{(}->}[r]  & {} } }

\newcommand{\dpp}{D_+^3}
\newcommand{\bpp}{B_+^3}
\newcommand{\dppd}{D_+^2(0,\d)}
\newcommand{\od}{[0,\d]}
\newcommand{\mnk}{\MM\NN(K)}
\newcommand{\mnsk}{\MM\NN(S(K))}

\renewcommand{\mf}{Morse function}

\newcommand{\smk}{S^3\sm K}
\newcommand{\ssmk}{S^4\sm S(K)}

\begin{document}

\title
[On the Morse-Novikov number for 2-knots]
{On the Morse-Novikov number for 2-knots}
\author{Hisaaki Endo  and  Andrei Pajitnov}
\address{Tokyo Institute of Technology,
Tokyo, Japan}
\email{}
\address{Laboratoire Math\'ematiques Jean Leray 
UMR 6629,
Universit\'e de Nantes,
Facult\'e des Sciences,
2, rue de la Houssini\`ere,
44072, Nantes, Cedex}                    
\email{andrei.pajitnov@univ-nantes.fr}

\thanks{} 
\begin{abstract}
Let $K\sbs S^4$ be a 2-knot, that is,
a smoothly embedded 2-sphere in $S^4$.
The Morse-Novikov number $\MM\NN(K)$ is 
the minimal possible number of critical points of a Morse map 
$S^4\sm K\to S^1$ belonging to the canonical class 
in $H^1(S^4\sm K)$. 
We prove that for a classical knot $K\sbs S^3$ 
the Morse-Novikov number of the spun knot $S(K)$
is $\leq 2\MM\NN(K)$. This enables  us to compute
$\MM\NN(S(K))$ 
for every classical knot $K$ with tunnel number 1. 
\end{abstract}
\maketitle
\tableofcontents

\section{Introduction}
\label{s:intro}
\subsection{Overview of the article}
\label{su:overv}

Let $K\sbs S^4$ be a 2-knot, that is, a $\smo$ embedding of $S^2$ into $S^4$.
We say that $K$ is fibred, if the complement $C_K=S^4\sm K$ admits a fibration over 
$S^1$, which is standard nearby $K$ 
(see Definition \ref{d:def-reg}).
In general a Morse map $C_K\to S^1$ has critical points, 
the minimal number of these critical points
will be called {\it the Morse-Novikov number of $K$}
and denoted $\MM\NN(K)$.
The aim of this paper is to study this invariant of 2-knots 
and to compute it for 
several families of knots,

The Novikov homology provides lower bounds for the number $\MM\NN(K)$, see 
Section \ref{s:lobounds}.
In Section \ref{s:saddle}
we introduce {\it the saddle number $sd(K)$} of a 2-knot $K$;
it is defined as a minimal possible number of critical points of index
1 of a generic projection of $K$ 
to a line in $\rr^4$. This number can be considered as a
2-dimensional analog of the bridge number of a classical knot.
We prove that $\mnk\leq 2sd(K)$.
Using the results of Section \ref{s:lobounds}
we deduce the following homological lower bound for 
the saddle number in terms of the Novikov torsion numbers:
$$
\wh q_1(C_K)+
\wh q_2(C_K)
\leq sd(K).
$$

In Sections \ref{s:spun} -- \ref{s:mn-spun} we study the Morse-Novikov
numbers of spun knots. We prove in particular
that for a classical knot $K\in S^3$ the Morse-Novikov number of the spun knot 
$S(K)$ satisfies
$$
\mnsk\leq 2\mnk.
$$

In \cite{P-t} the second author proved that if $K\sbs S^3$ is a classical knot, then 
$\mnk \leq 2t(K)$
where $t(K)$ is the tunnel number of $K$. Using this inequality we prove that 
$\mnsk \leq 4 $ for any classical knot $K$  of tunnel nuber 1 (Section \ref{s:mn-spun}).

The case of high-dimensional knots  is different from that of 1- and 2-knots. 
We gathered some results on $\mnk$ of knots of dimension $\geq 6$ in Section 
\ref{s:fibr-high-dim}. 
In particular we prove that for a given knot $K^n, n\geq 6$ 
and two numbers $p,q$ of the same parity
the fibredness of $S_p(K)$ is equivalent to the fibredness of $S_q(K)$.

\subsection{Definition and first properties of 
the Morse-Novikov numbers for 2-knots}
\label{s:defs}

Let $K\sbs S^4$ be the image of the $\smo$ embedding 
of $S^2$ to $S^4$. Choose a $\smo$ trivialisation 
$$
\Phi:N(K) \to  K \times B^2(0,\e)
$$
of a tubular \nei~ $N(K)$. 
We denote the complement to $K$ by $C_K$.
This space is non-compact, and to develop the 
Morse theory on it, we will assume that the 
functions and vector field have standard behaviour 
nearby $K$. For the Morse maps and their gradients we use the 
terminology of \cite{CVMT}.

\bede\lb{d:def-reg}
A Morse map $f ~:~C_K\to S^1$ is said to be {\it regular\/}
if there is  a $\smo$ trivialisation   
\begin{equation}\lb{f:triv}
\Phi:N(K) \to K \times B^2(0,\e)
\end{equation}
of a tubular neighborhood of $K$
such that the 
restriction $f|\:N(K)$
satisfies 
$f\circ \Phi^{-1} (x,z)=  z/|z|.$

An $f$-gradient $v$ of a regular Morse map $f~:~C_K\to S^1$
will be called {\it regular} if there is a $\smo$ trivialisation   
\rrf{f:triv} 
such that $\Phi^*(v)$ equals $(0,v_0)$
where $v_0$ is the Riemannian gradient of 
the function $z\mapsto  z/|z|$. 

A pair $(f,v)$, where $f$ is a regular Morse map and $v$ is its regular
gradient will be called {\it a Morse pair}.
\end{defi}

If $f$ is a Morse map of a manifold to 
$\RRR$ or to $S^1$,
then we denote by $m_p(f)$
the number of critical points of $f$ of index $p$.
The number of all critical points of $f$
is denoted by $m(f)$.

\bede\lb{d:def-MN}
The minimal possible number 
$m(f)$ where $f:C_K\to S^1$ is a 
regular Morse map is called 
{\it the Morse-Novikov number of $K$.}

A regular Morse map $f :C_K\to S^1$ is called {\it minimal\/}
if the number $m(f)$
is minimal on the class of all regular Morse functions. \end{defi}

The following proposition allows to get rid of 
the local maxima and minima for regular functions. 
Its proof 
repeats the proof of Lemma 3.2 in \cite{PRW} 
and will be omitted. 
\bepr\lb{p:ind0-4}
Let $f:C_K\to S^1$ be a 
regular Morse map. 
Then there is a regular Morse map
$g:C_K\to S^1$
such that 
$m_i(g)\leq m_i(f)$ for every $i$
and 
$m_0(g)=m_4(g)=0$.
\enpr

The next Corollary is immediate.
\beco\lb{c:no-min-max}
There is a minimal Morse function 
$f:C_K\to S^1$
without critical points of indices $0$ and $4$.
\enco

\bede\lb{d:strong-min}
A Morse map 
$f:C_K\to S^1$
is called strongly minimal 
if for every $i, 0\leq i \leq 4$ 
the number $m_i(f)$ is minimal 
in the class of all regular
functions $C_K\to S^1$.
\end{defi}

\bere\lb{r:min-str-min}
For all the 2-knots where we are able to compute the 
Morse-Novikov number,
the strongly minimal Morse functions exist
(see Section \ref{s:mn-spun}).
However 
it is not clear whether strongly minimal functions exist for every $K$.
The case of 2-knots is different here from the case of 
classical knots, for which the concept of the strongly 
minimal Morse map is the same as that of minimal map 
(by \cite{PRW}). 

\end{rema}

\bere\lb{r:multidim}
The main object of study in
the present paper are 2-knots. However the definitions 
above generalize immediately 
to the knots of any dimension, and we will use 
corresponding terminology
throughout the paper.
\enre

\section{Lower bounds from Novikov homology}
\lb{s:lobounds}

Let $L=\zz[t^{\pm }]$;  denote by  $\wh L=\zz((t))$ and 
$\wh L_\qq=\qq((t))$
the rings of all 
series in one variable $t$ 
with integer (resp. rational) coefficients  and finite negative part. 
Recall that $\wh L$ is a PID, and $\wh L_\qq$ is a field.

Consider the infinite cyclic covering 
$\ove {C_K} \to C_K$; the Novikov homology of $C_K$ is defined as follows: 
$$
\wh H_*(C_K)= H_*(\ove {C_K})\tens{L}\wh L.
$$
The rank and torsion number of the $\wh L$-modle 
$\wh H_*(C_K)$
will be denoted by $\wh b_k(C_K)$, resp. $\wh q_k(C_K)$.
For any regular Morse function
$f:C_K\to S^1$ there is a Novikov complex 
$\NN_*(f)$ 
 over $\wh L$ generated in degree $k$ by critical points of $f$ of index $k$
 and such that 
 $H_*(\NN_*(f))\approx \wh H_*(C_K).$
 We deduce the Novikov inequalities
 $$  
 \sum_k \Big(\wh b_k(C_K)+\wh q_k(C_K)+\wh {q}_{k-1}(C_K)\Big)
 \leq
  \mnk.
 $$
The numbers $\wh b_k(C_K), \ \wh q_k(C_K)$
satisfy certain relations. The homology of $C_K$ is the same as that of $S^1$, 
therefore (by \cite{Mi}) the $\qq$-vector space 
$H_*(\ove{C_K}, \qq)$ has finite dimension. This implies 
$H_*(\ove{C_K})\tens{L}\wh L_\qq=0$, and $\wh b_i(C_K)=0$ for all $i$.
It is clear that $\wh q_0=\wh q_4=0$.
Furthermore, since there is always a regular Morse map 
without critical points of index $4$,
the Novikov homology has no torsion in degree $3$, therefore $\wh q_3=0$.
Thus the Novikov inequalities boil down to the following:
\bq\lb{f:mn-q}
2\Big(\wh q_1(C_K)+\wh q_2(C_K)\Big)  \leq \mnk.
\end{equation}
Observe that both $\wh q_1, \wh q_2$
can be non-zero, as the example of the spun knot of the $5_2$-knot shows. 

In the sequel we will use the universal Novikov complex
as well.
Denote by $G$ the group $\pi_1(C_K)$, put 
$\L=\zz G$.  
Let $\xi\in H^1(M,\zz)$ be the  generator of 
the group $H^1(M,\zz)\approx \zz$,
positive on every meridian of $K$, we will 
call it {\it the canonical generator}. It can be 
considered as a \ho~ $G\to\zz$.
Recall the Novikov ring
$$
\Lxi=
\{\l=\sum_{k\in \nn} n_k g_k ~~|~~ n_k\in\zz, 
g_k\in G\ {\rm and }\ \xi(g_k)\to -\infty \}.
$$
For every \mf~ $f:C_K\to S^1$ we have a chain complex 
$\wi\NN_*(f)$ of free $\Lxi$-modules generated in 
degree $k$ by $Crit_k(f)$, and such that 
$$
C_*(\wi C_K) \tens{\L}\Lxi \sim      \wi\NN_*(f)
$$
(see \cite{P-o}).
Tis chain complex is defined via counting flow lines of an $f$-gradient $v$, 
so we will denote it by $\wi\NN_*(f,v)$ when the dependance on $v$ is important
(as for example in  Subsection \ref{su:nov-com}).

\section{Motion pictures and saddle numbers}
\lb{s:saddle}

Let $K\sbs \rr^4$ be a 2-knot. 
Choose a  projection $p$
of $\rr^4$ onto a line. 
Assume that the critical points of the 
function $p|K$ are non-degenerate.

\bede\lb{d:saddle}
The minimal number of saddle points of the function $p|K:K\to\rr$
(the minimum is taken over all embeddings of $S^2$ in $\rr^4$
ambient isotopic to $K$), will be called 
{\it saddle number of $K$}, and denoted by $sd(K)$.
\end{defi}

It is not difficult to prove that for any classical knot $K$
we have
$$
sd(S(K)) \leq 2(b(K)-1),
$$
where $b(K)$ denotes the bridge number of $K$, and $S(K)$ is the spun knot of $K$.
The invariant $sd(K)$ is closely related to the {\it $ch$-index }
of $K$, introduced and studied by K. Yoshikawa in \cite{Yoshikawa1994}.
In particular, we have $sd(K)\leq ch(K)$. 
In order to relate the  number $sd(K)$ to $\mnk$ we 
will reformulate the definition of the saddle number.

%2.1 \ \  Motion picture and $2$-knot exterior
Let $K\subset S^4$ be a $2$-knot. 
The equatorial $3$-sphere $\Sigma^3$ of the standard Euclidean sphere $S^4$ 
divides $S^4$ into two parts: 
$$
S^4=D^4_+\cup D^4_-,\ \  {\rm with }\ \  D^4_+\cap D^4_-=\Sigma^3.
$$
We assume that $K$ is included in ${\rm Int}(D^4_-)$ and 
$K$ does not include the center of $D^4_-$. 
Perturbing the embedding $K\subset D^4_-$ if necessary, 
we can assume that the restriction $\rho=r|_K$ 
of the radius function $r:D^4_-\rightarrow [0,1]$ 
is a Morse function. 
The family $\{(r^{-1}(t), \rho^{-1}(t))\}_{t\in [0,1]}$ of possibly singular knots can be drawn 
as a {\it motion picture} (see \cite{Kamada2002}, Chapter 8). 
Each singularity of a knot in the family corresponds to a critical point of $\rho$. 
A critical point of $\rho$ of index $0$ ($1$, $2$, respectively) is called 
{\it minimal point} ({\it saddle point}, {\it maximal point}, respectively) of $\rho$, 
which is represented by a {\it minimal band} ({\it saddle band}, {\it maximal band}, respectively) 
in (a modification of) the motion picture. 

It is clear that 
the minimal number of the saddle points for all such Morse functions $\rho$ 
is equal to $sd(K)$.

\bepr\lb{p:mn-s} \ \ 
$\mathcal{MN}(K)\leq 2\, sd(K)$. 
\enpr
{\it Proof}. 
Since $\rho$ is a Morse function, 
the manifold $D^4_-\setminus {\rm Int}\, N(K)$ admits a handle decomposition 
with $m_i(\rho)$ $(i+1)$-handles for $i\in\{0,1,2\}$ (see \cite{GS1999}, Proposition 6.2.1). 
The exterior $E(K)=S^4\setminus {\rm Int}\, N(K)$ of $K$ is obtained by attaching 
a $4$-handle $D^4_+$ to $D^4_-\setminus {\rm Int}\, N(K)$. 
Since $D^4_-\setminus {\rm Int}\, N(K)$ is connected, 
there is a $3$-handle in $D^4_-\setminus {\rm Int}\, N(K)$ 
which connects $\partial N(K)$ with $\partial D^4_-$. 
Thus the $3$-handle cancels the $4$-handle $D^4_+$ (cf. \cite{Milnor1965}, Section 5). 
Turning the handlebody upside down, 
we obtain a dual decomposition of $E(K)$ and 
a corresponding Morse function $f:E(K)\rightarrow \mathbb{R}$ 
which is constant on $\partial E(K)$ and the following Morse numbers: 
$m_1(f)=m_2(\rho)-1$, $m_2(f)=m_1(\rho)$, $m_3(f)=m_0(\rho)$, $m_4(f)=1$. 

Using the argument from \cite{P-t}, p. 629,
we can deform the real-valued Morse function $f$
to a circle-valued regular function $\phi:E(K)\to S^1$,
such that $m_k(f)=m_k(\phi)$ for every $k$.
Consider the function $-\phi$, which has one critical point 
of index $0$. Applying the cancellation of this local minimum,
we obtain a Morse function $\psi:E(K)\to S^1$
belonging to the class $-\xi$, and such that such that 
$m_0(\psi)=0,\ m_1(\psi)=m_3(f)-1, \ m_2(\psi)=m_2(f), \ m_3(\psi)=m_1(f), \ m_4(\psi)=0$.
Put $g=-\psi$. Then we have
$$
m_0(g)=m_4(g)=0, \ \ 
m_1(g) = m_2(\r) -1, \ $$
$$
m_2(g) = m_1(\r), \ \ 
m_3(g) = m_0(\r) -1.
$$
Observe that $m_0(\r)-m_1(\r)+m_2(\r)=\chi(S^2)=2$, 
therefore the total number of critical points of 
$g$ equals $2m_1(\r)$. Choosing the function $\r$ with 
$m_1(\r)=sd(K)$ we accomplish the proof. $\qs$

Taking into account the inequality 
\rrf{f:mn-q}
we obtain the following.
\beco\lb{c:q-saddle}
$$
\wh q_1(C_K)+\wh q_2(C_K)
\leq sd(K)\leq ch(K). \ \ \ \qs
$$
\enco

\section{Circle-valued Morse maps for spun knots}
\lb{s:spun}

Let $K\sbs S^3$ be a classical knot, and $(\phi, v)$ be a 
Morse pair on $C_K$. Denote by $S(K)$ the spun knot of $K$
(see Subsection \ref{su:spun-def}
for definition). In this section we associate to $(\phi,v)$ 
a regular Morse pair $(F,w)$ on $C_{S(K)}$. We compute the Novikov 
complex of $(F,w)$
in terms of the Novikov 
complex of $(\phi,v)$
(Propositions \ref{p:F-Morse} and \ref{p:nov-com}).

\subsection{Spun knots: the definition}
\lb{su:spun-def}

Let us recall the classical construction of the 
spun knot for a classical knot $K\sbs S^3$ (due to Artin).
The equatorial 2-sphere $\Sigma$ of the standard Euclidean 
sphere $S^3$ divides 
$S^3$ into two parts:
$$S^3=D^3_+\cup D^3_-,\ \  {\rm with }\ \  D^3_+\cap D^3_-=\Si.$$
We can assume that $K\cap D^3_-$ is a half-circle of the sphere $S^3$.
Let 
$K_+=D^3_+\cap K$, observe that 
$\Int(D^3_+\sm K_+)$
is diffeomorphic to $S^3\sm K$.
The sphere $S^4$ can be considered as an open book with binding 
$\Si=\pr \dpp$
and pages diffeomorphic to $\dpp$.
In other words, $S^4$ is obtained
by rotating of $\dpp$ around $\pr\dpp=\Si$. 
The result of rotation of $K_+\sbs \dpp$ is an embedded 
2-sphere in $S^4$,
thus a 2-knot, which is called {\it the spun knot of $K$};
we denote it by $S(K)$.

\subsection{Geometric setup}\lb{su:gs}
${}$
\pa
{\it 1. The embedding $K_+\sbs \dpp$.}

The intersection $K_+\cap\Si$ consists of two points,
denote them by $s$ and $n$. 
We can assume that there is a collar
$h:\Si\times [0,\d]\to \dpp$
such that 
$h^{-1}(K_+\cap \dpp)=\{s,n\}\times [0,\d]$.
Let $(\phi,v)$ be a Morse pair on $C_K$.
Without changing the Novikov complex $\NN_*(\phi,v)$ we can 
make the following assumptions on the map $\phi$ and its gradient $v$:
\been\item
All the critical points of $\phi$ are in $\dpp$.
\item 
$\phi\circ h$ is constant  along each interval
$x\times [0,\d]$, where $x\in\Si$.
\item There is a local coordinate system 
$Q_n:\dppd\to\Si$ around $n$ such that 
$(h\circ (Q_n\times \Id)^{-1})^*(v)$ equals the vector field 
$(\frac{z}{|z|}, 0)$ in  $\dppd\times [0,\d]$.
There is a similar coordinate system $Q_s$ around $s$.
\item
The trivialisation 
$
N(K_+)
\arrr \r 
\dppd\times K_+
$ required by the definition 
of regular Morse pair 
is compatible with $h$
in the \nei~ 
$h(\dppd\times [0,\d])$.
\enen

In the sequel  we will need two auxiliary 
functions defined on $\dpp$.
Let $\xi:[0,\d]\to\rr$ be a $\smo$ function
such that $\xi(r)=0$ for $r\in [0,\d/3]$, and $\xi'(r)>0$ for
$r\in ]\d/3, 2\d/3[$, and $\xi(r)=1$ for
$r\in [2\d/3, \d]$.

A) \ \ Define a function $\a:\dpp\to [0,1]$ as follows.
For $x=h(y,r)$ with $y\in 
\Si\times [0,\d]$
put $\a(x)=\xi(r)$. For $x$ outside $\Im h$ put $\a(x)=1$.
Then $\a$ is a $\smo$ function vanishing in a \nei~ of $\Si$.

B) \ \ Define a function $\b:\dpp\to [0,1]$ as follows.
For $x=\r(z,t)$ with
$z\in \dppd, t\in K_+$
put $\b(x)=\xi(|z|)$. Otherwise put $\b(x)=1$.
Then $\b$ is a $\smo$ function vanishing in a \nei~ of $K$.

${}$
\pa
{\it 2. The embedding $S(K)\sbs S^4$.}

The \nei~ of $\Si$ in $S^4$ can be parameterized by a diffeomorphism
$H:\Si\times \od\times S^1\to N(\Si)$ where
$H~|~\Si\times \od\times \{1\}=h$
and for a given $a\in H$ the map 
$(a,r,\t)\mapsto H(a,r,\t)$ gives the polar 
coordinates in the 2-disc normal to $H(a,0,0)\in\Si$.
The coordinate $\t$ is defined as a $\smo$ 
function on $S^4\sm \Si$. We have a diffeomorphism
$$S^4\sm \Si\approx \bpp\times S^1,$$
its  second projection
$S^4\sm \Si \to S^1$ extends the angle cordinate $\theta$ defined on 
$N(\Si)\sm \Si$, its first projection
$S^4\sm \Si\approx \bpp\times S^1\to \bpp$
will be denoted by $p_1$.

\subsection{Construction of a Morse map $F:S^4\sm S(K)\to S^1$}
\lb{su:morse-map}

Define a function $F_0$ on $S^4\sm S(K)$ by
$F_0(u,\t)=\phi(u)$. In other words, $F_0=(\phi\circ p_1)(x)$,
where $p_1:S^4\sm \Si\to \bpp$ is the projection of the first 
factor of the cartesian product.
It is clear that $F_0$ is $\smo$ on $S^4\sm (S(K)\cup \Si)$.
Using our assumptions on $\phi$ 
(see Subsection \ref{su:gs})
we deduce that $F_0$ extends to a $\smo$ function on the whole 
of $S^4\sm S(K)$. We have 
${\rm Crit } F_0\approx S^1\times {\rm Crit } \phi$,
therefore $F_0$ is not a Morse function.
We will now construct a small perturbation of $F_0$
which will have only non-degenerate critical points.
Let $h:S^1\to\rr$ be a Morse function with two 
non-degenerate critical points
(e.g. $h(\t)=\sin(\t)$). 
Define a function
$\wi h: S^4\sm \Si\to S^1$ by 
$\wi h(\t, u)=h(\t)$.

Extend the functions $\a, \b: \dpp\to[0,1]$
constructed in Subsection \ref{su:gs} 1)
to the functions on $S^4$ invariant 
\wrt~ action of $S^1$. 
We will denote the resulting functions by the same 
symbols $\a,\b$. 
Then $\a$ is a $\smo$ function, 
vanishing in a \nei~ of $\Si$ and equal to $1$
outside $N(\Si)$.
The function   $\b$ vanishes in a \nei~ of 
$S(K)$ and equals
$1$ outside $N(S(K))$. Define a function $G:S^4\to S^1$ as follows:
\begin{align*}
G(x)=\a(x)\b(x)\wi h(x)  &\ \ {\rm for } \ \ x\notin \Si; \\
G(x)=0  &\ \ {\rm for } \ \ x\in \Si. 
\end{align*}
Then $G$ is a $\smo$ function vanishing in a \nei~ of 
$\Si\cup S(K)$ and equal to $1$
outside $N(\Si)\cup N(S(K))$.

\bepr\lb{p:F-Morse}
For any $\e$ sufficiently small, 
the function $F$ defined by
$$F(x) =F_0(x)+\e G(x), $$
is a Morse function, and 
$$m_i(F)= m_i(\phi)+m_{i-1}(\phi)
\ \ {\rm for\ every } \ \ i.$$
\enpr
\Prf
In the domain 
$N(\Si)\cup (N(S(K)) \sm S(K))$
the norm of $F_0'$ is bounded from below
by a strictly positive constant,
therefore the gradient of $F$ is non-zero everywhere 
in this domain if $\e$ is sufficiently small.
In the domain
$$
S^4 \sm 
\Big(N(\Si)\cup N(S(K))\Big) 
\approx (\bpp\sm K_+)\times S^1
$$
the function $F$ is 
diffeomorphic to the function 
$\wi F:(\bpp\sm K_+)\times S^1\to S^1$
defined by
$$
\wi F(u,\t)=\phi(u)+\e h(\t).
 $$            
 The proposition follows. $\qs$

\subsection{The Novikov complex of $F$}
\lb{su:nov-com}

We will now construct a suitable gradient for $F$.
Let $v$ be a $\phi$-gradient on $\bpp\sm K$ satisfying the assumptions
from Subsection \ref{su:gs}.
Define a vector field $\wi v$ on the product 
$(\bpp\sm K)\times S^1$ by
$\wi v(x,\t)=v(x)$, where $\t\in S^1, x\in \bpp$.
Carry over the vector field $\wi v$ to $S^4\sm S(K)$
(we will keep the same notation for the resulting field).
Let $u$ be any gradient for $h:S^1\to\rr$.
Define a vector field $\wi u$ on the product 
$\bpp\times S^1$ by the formula
$\wi u(z,t)=u(t)$
(where $z\in \bpp, t\in S^1$).
Carry over this vector field 
to $S^4\sm \Si$ and define a vector field $\bar u$
on $S^4\sm \Si$ setting
$\bar u(x)=\a(x)\b(x)\wi u(x)$. 
The vector field  $\bar u$ extends to a $\smo$ vector field on $S^4$
(which will be denoted by the same symbol $\bar u$),
vanishing in a \nei~ of $\Si\cup S(K)$. 
It is clear that for $\e$ small enough,
the vector field $w= \bar v + \e \bar u$
is an $F$-gradient.

\bede\lb{d:susp}
For a chain complex $A_*$ over a ring $R$ we denote
by $\s A_*$ the suspension of $A_*$,
that is, $(\s A_*)_k= A_{k-1}$.
\end{defi}
In the next proposition we compute the Novikov complex of $F$
in terms of the Novikov complex of $\phi$. Observe that
the base ring of both complexes 
is the same, it is the Novikov completion of $\zz G$, where $G=\pi_1(\smk)
\approx \pi_1(\ssmk)$.

\bepr\lb{p:nov-com}
We have
$$\NN_*(F, w)\approx 
\NN_*(\phi, v)
\oplus
\s \NN_*(\phi, v).$$
\enpr
\Prf
It is easy to check 
(using our assumptions on $v$)
that any flow line of $w$ joining critical points of $F$ does not intersect 
the subset 
$N(K)\cup N(\Si)$, and therefore remains in the domain,
where the vector
field $w$ is the direct product of vector fields 
$\e \wi u (x)$ and $\wi w(x)$. Therefore 
$$
\NN_*(F, w)\approx
\NN_*(\phi, v)
\otimes
\MM_*(h, u)
$$
(where $\MM_*(h,u)$ is the Morse complex for $h$)
and the Proposition follows. $\qs$

\subsection{Superspinning}
\lb{su:superspin}

The classical Artin construction has several generalizations,
in particular {\it the superspinning}.
This construction, due to
E.C. Zeeman
\cite{Z}
and D.B.A. Epstein
\cite{Ep}, 
associates to any $n$-knot
$K^n\sbs S^{n+2}$ and a natural number $p\geq 1$
an $(n+p)-$ knot in $S^{n+p+2}$
(see \cite{Fr}, p. 196). 
We will denote the resulting knot by $S_p(K)$.
The results of the previous subsection 
generalize directly to
the superspinning case. Let $\phi:S^{n+2}\sm K^n\to S^1$
be a regular Morse map, and $v$ a regular $\phi$-gradient.

\beth\lb{t:sup-spin-nov}
There is a regular Morse function
$F:S^{n+p+2}{\setminus} S_p(K)
\to S^1$
and a regular $F$-gradient $w$ such that
\been\item
$m_i(F)=m_i(\phi)+m_{i-1}(\phi)$,
\item
$\NN_*(F,w)=\NN_*(\phi,v)\oplus \s^p \NN_*(\phi,v)$.
\enen
\enth
\Prf The argument repeats the proof of Propositions 
\ref{p:F-Morse}
and \ref{p:nov-com}
with minor modifications. $\qs$

\section{Morse-Novikov numbers for spun knots}
\lb{s:mn-spun}

In this section we gathered some consequences 
of the constructions
developed in Section \ref{s:spun}.
The next Corollary is immediate from 
Proposition \ref{p:F-Morse}.

\beco\lb{c:ineq}
$\MM\NN(S_p(K))\leq 2\MM\NN(K)$. \ \ $\qs$
\enco
In particular, if $K$ is fibred, then $S_p(K)$ is fibred. 
The case of fibred knots was observed in 
\cite{AS}, together with the inverse implication:
if $K$ is a classical knot and 
$S(K)$ is fibred, then $K$ is fibred.
This last property is not valid in higher dimensions,
as an example of C. Kearton 
\cite{K} shows.

\beth\lb{t:mn2}
Let $K\sbs S^3$ be a classical knot with $\MM\NN(K)=2$.
Then $\MM\NN(S(K))=4$ and there is a strongly 
minimal Morse function on $S^4\sm S(K)$.
\enth
\Prf
Let $\phi:S^3\sm K\to S^1$ be  a regular Morse function
without local maxima or minima, 
with $m_1(\phi)=m_2(\phi)=1$. By Proposition \ref{p:F-Morse}
there is a 
regular Morse function
$\phi:S^4\sm S(K)\to S^1$
with
$m_1(F)=m_3(F)=1,\  m_2(F)=2, \ 
m_0(F)=m_4(F)=0$. Therefore 
$\mnsk\leq 4$. We will show that $F$ is 
actually a strongly minimal Morse map.

Let $H:S^4\sm S(K)\to S^1$ be any \mf. 
Denote the fundamental group of $S^3\sm K$ by $G$.
It is known that 
$\pi_1(\ssmk)\approx G$.
If $m_1(H)=0$, then a standard Morse-theoretic argument 
applied to the infinite cyclic cover $\ove{\smk}$ of $\smk$ 
implies that $\pi_1(\ove{\smk})$ 
is finitely generated, which is impossible, since $K$ is not fibred.
Therefore $m_1(H)\geq 1$.
A similar argument shows that $m_3(H)\geq 1$.
Assume now that $H$ is a minimal \mf.
Then 
$$m_0(H)=
m_4(H)=0, \ m_1(H)\geq 1, m_3(H)\geq 1,$$
and this implies $m_2(H)\geq 2$ 
(since $\chi(\ssmk)=0$).$\qs$

In the work \cite{G2} based on his earlier paper \cite{G1}
H. Goda proved that $\mnk=2$ for every non-fibred 
prime knot with $\leq 10$ crossings. 
Therefore $\mnsk=4$ for these knots.

In the work \cite{P-t} the second author 
proved that for any classical
knot $K$ we have 
\begin{equation}\lb{f:tun}
\mnk\leq 2t(K),
\end{equation}
where $t(K)$ is the tunnel number of $K$.
\beco\lb{c:tunn1}
Let $K\sbs S^3$ be a non-fibred classical knot with
tunnel number 1. Then $\mnsk=4$. $\qs$
\enco
\Prf 
By \rrf{f:tun}
we have 
$$\mnk\leq 2$$
on the other hand $S(K)$ is not fibred since $K$ is non-fibred. 
Hence $\mnsk=4$ by  Theorem \ref{t:mn2}. $\qs$

These results allow to compute the numbers
$\mnsk$ for many classical knots $K$. 
In the paper \cite{MSY} K. Morimoto, M. Sakuma and Y. Yokota
explicited many examples of tunnel number one knots, in particular,
an infinite series of Montesinos knots.
They prove that the Montesinos knot 
$K=M(3, (\a_1, \b_1),(\a_2, \b_2),(\a_3, \b_3))$
has the tunnel number 1, if $r=3$, and 
$\b_2/\a_2\cong \b_3/\a_3 \mod \zz$ and 
$b - \sum_{i=1}^{3}\b_i/\a_i =   1/(\a_1\a_2)$.

\beco\lb{c:monte}
Let $K$ be a non-fibered Montesinos knot satisfying the conditions above.
Then $\mnsk=4$. $\qs$
\enco
\bere\lb{r:HM}
In the  work \cite{HM} of M. Hirasawa and K. Murasugi 
it is shown that the fibredness of  most Montesinos knots with
tunnel number  1 
is detected by the monicness of the Alexander polynomial.
\enre

\section{On fibering of high-dimensional spun knots}
\lb{s:fibr-high-dim}

Let $K^n\sbs S^{n+2}$ be a knot with $n\geq 6$.
Let $C_K$ denote the complement of $K^n$ in $S^{n+2}$, 
and $\wi C_K$ be the universal covering of $C_K$.
Let $\xi\in H^1(C_K,\zz)\approx \Hom(\pi_1(C_K),\zz)$
be the canonical generator of the cohomology of $C_K$.
The knot $K^n$ is fibred if and only 
if the following two conditions hold
(see \cite{P-o}, \cite{L}):
\been
\item[F0)] 
The subgroup $\Ker\xi$ of $\pi_1(C_K)$ is finitely presented.
\item[F1)]
The Novikov homology $\wh H_*(\wi C_K)$ vanishes,
\item[F2)]
The Whitehed torsion $\tau(K)\in \Wh(\Lxi)$ of the completed chain complex
$$
C_*(\wi C_K)\tens{\L} \Lxi
$$
is equal to $0$.
\enen 
\noindent
The results obtained in the previous sections 
allow to compare the conditions F0) -- F2) for
a knot $K$ and its spun knots. For $m\geq 1, p\geq 1$
we denote by  $S^m_p(K)$ the result of $m$ 
iterations of the $p$-spinning construction
of $K$.
Observe that the fundamental groups of the complement to the knot 
is isomorphic to that of the complement to its spun knot,
therefore the condition F0) holds for $K$ if and 
only if it holds for $S^m_p(K)$ (with any $p, m$).
Let $p\geq 1$. Theorem  \ref{t:sup-spin-nov}
implies that 
$$
\wh H_*(C_{S(K)}) \approx \wh H_*(C_{K})\oplus \s^p\wh H_*(C_{K})
$$
therefore the Novikov homology of $C_K$ vanishes if and only if the 
Novikov homology of $C_{S_p(K)}$ vanishes. 
The situation with the torsion is different, 
since Theorem  \ref{t:sup-spin-nov} 
implies
$$\tau(S_p(K))
=
(1+(-1)^p)\tau(K). $$
In particular, $\tau(S_p(K))=0$ for any $K$ if $p$ is odd, and 
$\tau(S_p(K))=2\tau(K)$ if $p$ is even. 

The next two propositions are now immediate.
\bepr\lb{p:iter-spin1}
If $p\geq 1$ and $p'\geq 1$ have the same parity, then 
$S_p(K)$ is fibred if and only if $S_{p'}(K)$ is fibred.
$\qs$
\enpr

Let $S^m_p(K)$ denote the result of $m$ 
iterations of the $p$-spinning construction
of $K$.

\bepr\lb{p:iter-spin2}
If $p$ is odd and $l, m\geq 1$ then 
$S^m_p(K)$
is fibred if and only if 
$S^l_p(K)$
is fibred. $\qs$
\enpr

\bere\lb{r:spun-tors}
If for some  knot $K^n\sbs S^{n+2}$
the conditions 
F0), F1)  above hold, but F2) does not hold,
then $S(K)$ is fibred, and $K$ is not fibred. 
However we do not know if such knots exist.
\enre

\section{Open questions}
\lb{s:open}

{\bf 1.} Is it true that for any 2-knot $K$ there exists a strongly minimal 
\mf~ $C_K\to S^1$?
(This is true for spun knots, see Theorem \ref{t:mn2}.)

{\bf 2.}  Is it true that for any classical knot $K$ we have 
$\mnsk=2\mnk$?
(This is true for any classical knot $K$ with $\mnk=2$, again by Theorem \ref{t:mn2}.)

{\bf 3.}  Is it true that for a knot $K$ of dimension $\geq 6$ we have $\mnk=\mnsk$?
In partiular is it true that $K$ is fibred if and only if $S(K)$ is fibred?

{\bf 4.}  It is not difficult to prove that 
for knots $K_1, K_2$ of any dimension
$\MM\NN(K_1\krest K_2)\leq \MM\NN(K_1)+ \MM\NN(K_2)$
(the argument  repeats the proof for the classical knots,
see \cite{PRW}). Is it true that 
$$
\MM\NN(K_1\krest K_2) =  \MM\NN(K_1)+ \MM\NN(K_2)
$$
for 2-knots? 

{\bf 5.} What is the relation between the saddle number and 
the unknotting number of a 2-knot? Is it true that $sd(K)\leq u(K)$?

\section{Acknowledgements}
\lb{s:acknow}

We began working on this project in 2014 during the visit of the first author to the 
Nantes University. The work was accomplished when the second author 
was visiting the Tokyo Institute of Technology in 2015.
The first author thanks the Nantes University and the 
GeanPyl program for the support and warm hospitality.
The second author thanks the Tokyo Institute of Technology
for support and  warm hospitality. 
The first author was partially supported by JSPS KAKENHI
Grant Number 25400082.

\end{document}